\documentclass[12pt,twoside]{article}

\usepackage{amssymb}
\usepackage{amsmath}
\usepackage{bbm}
\usepackage{mathrsfs}
\usepackage{xypic}

\makeatletter
\def\hsmash{\relax 
  \ifmmode\def\next{\mathpalette\mathhsm@sh}\else\let\next\makehsm@sh
  \fi\next}
\def\makehsm@sh#1{\setbox\z@\hbox{#1}\finhsm@sh}
\def\mathhsm@sh#1#2{\setbox\z@\hbox{$\m@th#1{#2}$}\finhsm@sh}
\def\finhsm@sh{\wd\z@\z@ \box\z@}
\makeatother

\makeatletter
\gdef\th@mychange{\normalfont\slshape
   \def\@begintheorem##1##2{\item
        [\hskip\labelsep \theorem@headerfont ##2. ##1  \,--\!--\!--\!--  ]}%
 \def\@opargbegintheorem##1##2##3{%
   \item[\hskip\labelsep \theorem@headerfont ##2. ##1\ {\upshape(}##3{\upshape)}. \,-----  ]}}
\makeatother

\RequirePackage{theorem}
\theoremstyle{mychange}

{\theorembodyfont{\rmfamily}\newtheorem{ttt}{}[section]}
{\theorembodyfont{\rmfamily}\newtheorem{ex}[ttt]{Example.}}
{\theorembodyfont{\rmfamily}\newtheorem{rem}[ttt]{Remark.}}
{\theorembodyfont{\rmfamily}\newtheorem{rems}[ttt]{Remarks.}}
{\theorembodyfont{\rmfamily}}

{\theorembodyfont{\rmfamily}\newtheorem{remo}[ttt]{Remark}}

{\theorembodyfont{\itshape}}
{\theorembodyfont{\itshape}\newtheorem{subl}[ttt]{Sublemma.}}
{\theorembodyfont{\itshape}\newtheorem{lem}[ttt]{Lemma.}}
{\theorembodyfont{\itshape}}
{\theorembodyfont{\itshape}\newtheorem{theo}[ttt]{Theorem.}}
{\theorembodyfont{\itshape}\newtheorem{coro}[ttt]{Corollary.}}

{\theorembodyfont{\rmfamily}\newtheorem{exo}[ttt]{Example}}

{\theorembodyfont{\itshape}\newtheorem{lemo}[ttt]{Lemma}}
{\theorembodyfont{\itshape}\newtheorem{propo}[ttt]{Proposition}}

\newcounter{abc}
\newenvironment{abc}{\begin{list}{\rm \alph{abc}) }{\usecounter{abc} \leftmargin=0.0pt \labelsep=0.0pt \listparindent=0.0pt \labelwidth=0.0pt \parsep=\smallskipamount \itemsep=0.0pt \topsep=0.0pt \partopsep=\smallskipamount}}{\end{list}}
\newcounter{iii}
\newenvironment{iii}{\begin{list}{\rm \roman{iii}) }{\usecounter{iii} \leftmargin=0.0pt \labelsep=0.0pt \listparindent=0.0pt \labelwidth=0.0pt \parsep=\smallskipamount \itemsep=0.0pt \topsep=0.0pt \partopsep=\smallskipamount}}{\end{list}}

\newcommand{\calA}{\mathscr{A}}
\newcommand{\calO}{\mathscr{O}}
\newcommand{\calS}{\mathscr{S}}

\newcommand{\calHC}{\mathscr{HC}}

\newcommand{\frakf}{\mathfrak{f}}
\newcommand{\frakl}{\mathfrak{l}}

\newcommand{\frakp}{\mathfrak{p}}
\newcommand{\frakt}{\mathfrak{t}}
\newcommand{\frakx}{\mathfrak{x}}

\newcommand{\frakL}{\mathfrak{L}}
\newcommand{\frakP}{\mathfrak{P}}

\newcommand{\bbA}{{\mathbbm A}}
\newcommand{\bbC}{{\mathbbm C}}
\newcommand{\bbF}{{\mathbbm F}}
\newcommand{\bbN}{{\mathbbm N}}
\newcommand{\bbQ}{{\mathbbm Q}}
\newcommand{\bbR}{{\mathbbm R}}
\newcommand{\bbZ}{{\mathbbm Z}}

\newcommand{\Pb}{\mathop{\text{\bf P}}\nolimits}

\newcommand{\HC}{\mathop{\text{\rm HC}}\nolimits}
\newcommand{\N}{\mathop{\text{\rm N}}\nolimits}
\newcommand{\Gal}{\mathop{\text{\rm Gal}}\nolimits}
\newcommand{\lcm}{\mathop{\text{\rm lcm}}\nolimits}

\newcommand{\C}{\mathop{\text{\rm C}}\nolimits}
\newcommand{\Cl}{\mathop{\text{\rm Cl}}\nolimits}

\newcommand{\reg}{\text{\rm reg}}

\renewcommand{\div}{\mathop{\text{\rm div}}\nolimits}

\renewcommand{\mod}{\text{ {\rm mod} }}

\newcommand{\ratarrow}{$%
$\definemorphism{rat}\dashed\tip\notip%
\spreaddiagramcolumns{-12pt}%
 - \!\!\diagram%
\rrat & 
\enddiagram\!\!$%
$}

\newcommand{\br}{ }
\newcommand{\brr}{, }

\def\rightend#1#2{{%
 \leavevmode\nobreak\hskip .5em plus 1fil
 \penalty600 \hskip 0pt plus -1filll
 \vadjust{}\nobreak\hskip 0pt plus 1filll%
 #1\parfillskip=#2\relax \par}}

\def\eop{\ifmmode\rule[-22pt]{0pt}{1pt}\ifinner\tag*{$\square$}\else\eqno{\square}\fi\else\rightend{$\square$}{0pt}\fi}

\renewcommand{\thefootnote}{\arabic{footnote}}
\author{Andreas-Stephan Elsenhans${}^*$ and J\"org Jahnel${}^\ddagger$}

\date{}

\title{Cubic surfaces violating the Hasse principle \\are Zariski dense in the moduli scheme}

\begin{document}
\renewcommand{\thefootnote}{\fnsymbol{footnote}}

\maketitle

\begin{abstract}
We construct new examples of cubic surfaces, for which the Hasse principle fails. Thereby, we show that, over every number field, the counterexamples to the Hasse principle are Zariski dense in the moduli scheme of non-singular cubic surfaces.
\end{abstract}

\footnotetext[1]{School of Mathematics and Statistics F07, University of Sydney, NSW 2006, Sydney, Australia,\\
{\tt stephan@maths.usyd.edu.au}, Website: {\tt http://www.staff.uni-bayreuth.de/$\sim$bt270951}}

\footnotetext[3]{\mbox{D\'epartement Mathematik, Universit\"at Siegen, Walter-Flex-Str.~3, D-57068 Siegen, Germany,} \\
{\tt jahnel@mathematik.uni-siegen.de}, Website: {\tt http://www.uni-math.gwdg.de/jahnel}}



\section{Introduction}

\begin{ttt}
Cubic surfaces
over~$\bbQ$
that violate the Hasse principle are known for more than 50~years. The~first example of a cubic surface, for which the Hasse principle provably fails, was contrived by Sir~Peter Swinnerton-Dyer~\cite{SD62}, in~1962. The~construction had soon been generalized by L.\,J.~Mordell~\cite{Mordell}, who found a whole family of~examples. A~further generalization was recently given by one of the authors~\cite{J}.

A~completely different kind of counterexample, being a diagonal cubic surface with a very particular coefficient vector, was discovered by J.\,W.\,S.~Cassels and M.\,J.\,T.\ Guy~\cite{CG}, in~1966. Later,~J.-L.~Colliot-Th\'el\`ene, D.~Ka\-nev\-sky, and \mbox{J.-J.}~Sansuc~\cite{CTKS} studied the arithmetic of these surfaces, in~general.

Somewhat~surprisingly, it seems that, until today, no cubic surface has been found that violates the Hasse principle without being of one of these two~types. On~the other hand, it is known that the Hasse principle is always valid for singular cubic surfaces~\cite{Skolem}.
\end{ttt}

\begin{ttt}
The~coarse moduli scheme of non-singular cubic surfaces over a base
field~$K$
is the complement of a hypersurface in the four-dimensional weighted projective space
$\Pb(1,2,3,4,5)_K$
\cite{Clebsch}, \cite[Section~9.4.5]{Dolgachev}. All~diagonal cubic surfaces are geometrically isomorphic to each~other. Thus,~they correspond to a single point on the moduli~scheme.

On the other hand, the Swinnerton-Dyer-Mordell type surfaces are contained in a two-dimensional closed subscheme
of~$\Pb(1,2,3,4,5)_\bbQ$.
Indeed,~they are given by equations of the~form
\begin{equation}
\label{eins}
T_3(a_1T_0 + d_1T_3)(a_2T_0 + d_2T_3) = \N_{K/\bbQ}(T_0 + \theta T_1 + \theta^2 T_2) \, ,
\end{equation}
for
$K/\bbQ$
a cyclic cubic field extension,
$\theta \in K$,
and
$\N_{K/\bbQ}$
the norm map from
$K$
to~$\bbQ$.

Such~surfaces have at least three Eckardt~points. The~reason is that the three tritangent planes
$V(T_3)$,
$V(a_1T_0 + d_1T_3)$,
and~$V(a_2T_0 + d_2T_3)$
have a line in~common. Thus,~on each of the three tritangent planes
$V(T_0 + \theta^{(i)} T_1 + \theta^{(i)} {}^2 T_2)$,
the corresponding three lines meet at a single~point. Lemma~\ref{dim2} implies the~claim.
\end{ttt}

\begin{rems}
\begin{iii}
\item
Calculations~with concrete coefficients indicate that, gener\-i\-cally, there are not more than three Eckardt points on the surfaces~(\ref{eins}).

In~this case, the automorphism groups of the cubic surfaces are isomorphic
to~$S_3$
\cite[Proposition~9.1.27]{Dolgachev}. Therefore,~the surfaces are of type~VIII in I.\,V.~Dolgachev's classification~\cite[Table~9.6]{Dolgachev}.
\item
On~the other hand, a diagonal cubic surface has 18 Eckardt points, which is the maximal number a non-singular cubic surface may~have.
\end{iii}
\end{rems}

\begin{ttt}
Let
$\calHC_K \subset \Pb^{19}(K)$
be the set of all cubic surfaces violating the Hasse principle and
$\C\colon \Pb^{19} \ratarrow \Pb(1,\ldots,5)$
the Clebsch invariant~map. In~view of the considerations above, one is tempted to consider the following~problems.

\begin{iii}
\item
Describe~the Zariski closure in the moduli space of the locus
$\HC_K := \C(\calHC_K)$
of the counterexamples to the Hasse~principle.
\item
If
$\dim \overline\HC_K < 4$
then find the geometric properties of cubic surfaces that are implied by the arithmetic property of being a counterexample to the Hasse~principle.

In~particular, does every cubic surface that does not fulfill the Hasse principle automatically have Eckardt~points?
\end{iii}
\end{ttt}

\begin{ttt}
In~this article, we will show that actually
$\overline\HC_K$
is the full moduli space. I.e.,~that the Hasse counterexamples are Zariski dense in the moduli space of cubic~surfaces. In~particular, Problem~ii) is~pointless. Although~certainly the case of the base
field~$\bbQ$
is of particular interest, we will work over an arbitrary number
field~$K$.
\end{ttt}

\section{A family of cubic surfaces}

\begin{ttt}
We~consider the cubic surface
$S$
over a number
field~$K$,
given by the equation
\begin{equation}
\label{zwei}
T_0T_1T_2 = \N_{L/K}(aT_0 + bT_1 + cT_2 + dT_3) \, ,
\end{equation}
for
$L/K$
a cyclic cubic field extension
and~$a,b,c,d \in L$.
\end{ttt}

\begin{ttt}
Let~us assume
$d \neq 0$
as, otherwise, this surface is a~cone. Further,~we require that
$a/d$,
$b/d$,
$c/d \not\in K$.
Then, as
$\N_{L/K}$
represents zero only trivially,
$S$
has no
$K$-rational
point
$(t_0 \!:\! t_1 \!:\! t_2 \!:\! t_3) \in S(K)$
with more than one of
$t_0, t_1, t_2$
being equal
to~$0$.

To~prove
$S(K) = \emptyset$
for particular choices of
$K$,
$a$,
$b$,
$c$,
and~$d$,
our strategy is as~follows. Suppose~that there is a point
$(t_0 \!:\! t_1 \!:\! t_2 \!:\! t_3) \in S(K)$.
Among~$t_1/t_0$,
$t_2/t_1$,
and
$t_0/t_2$,
consider an
expression~$q$
that is properly defined and non-zero.

We~will show that, for every prime ideal
$\frakl$
of~$K$
with the exception of exactly one,
$q \in K_\frakl$
is in the image of the norm map
$N \colon L_\frakL \to K_\frakl$,
for~$\frakL$
a prime
of~$L$
lying
above~$\frakl$.
Such~a behaviour, however, is incompatible with global class field theory, cf.~\cite[Chapter~VI, Corollary~5.7]{Neukirch} or~\cite[Theorem~5.1 together with 6.3]{Ta}.
\end{ttt}

\begin{rem}
Equation~(\ref{zwei}) is similar to the Swinnerton-Dyer-Mor\-dell type. The~only difference is that the three linear forms
$T_0, T_1, T_2$
are linearly independent.
\end{rem}

\begin{ttt}
There is a conjecture due to J.-L.~Colliot-Th\'el\`ene~\cite[Conjecture~C]{CTS} that actually every cubic surface violating the Hasse principle does so via the Brauer-Manin obstruction, as introduced by Yu.\,I.~Manin in~\cite[Chapter~VI]{Ma}. The~new examples are in agreement with Colliot-Th\'el\`ene's~conjecture.

In~fact, the choice of a rational function such as
$T_1/T_0$
is not at all~arbitrary. The~principal divisor
$\div(T_1/T_0)$
is the norm of a divisor, which is the difference of two~lines. Thus,~the cyclic algebra
$$\textstyle \calA := L(S)\{Y\} / (Y^3 - \frac{T_1}{T_0})$$
over the function field
$K(S)$,
where
$Y t = \sigma(t) Y$
for
$t \in L(S)$
and a fixed generator
$\sigma \in \Gal(L/K)$,
extends to an Azumaya algebra over
$S$~\cite[Proposition~31.3]{Ma}.
This~shows that we work, indeed, with a particular case of the Brauer-Manin obstruction.
\end{ttt}

\begin{rems}
\label{Brauer}
\begin{iii}
\item
The~quotients
$\smash{\frac{T_1}{T_0} / \frac{T_2}{T_1} = \frac{T_1^3}{T_0T_1T_2}}$
and
$\smash{\frac{T_1}{T_0} / \frac{T_0}{T_2} = \frac{T_0T_1T_2}{T_0^3}}$
are norms of rational~functions. Thus,~the three expressions
$T_1/T_0$,
$T_2/T_1$,
$T_0/T_2$
actually define the same Brauer class.
\item
The non-singular cubic surfaces of the form~(\ref{zwei}) have a pair of Galois-invariant Steiner trihedra~\cite[Fact~4.2]{EJ1}. Therefore,~their Brauer groups are
\mbox{$3$-torsion}
of order
$3$
or~$9$.
A~general procedure to calculate the Brauer-Manin obstruction to the Hasse principle or weak approximation for such surfaces was described in~\cite{EJ2}.

The~methods developed there are, however, not necessary here in their full~strength. In~fact, the most complicated case treated in~\cite{EJ2} is that of an orbit structure of
type~$[9,9,9]$
on the 27~lines. The~surfaces of the form~(\ref{zwei}) generically have an orbit structure of type
$[3,3,3,9,9]$,
which is a technically much simpler~case.
\end{iii}
\end{rems}

\begin{ttt}
Although certainly Brauer classes work behind the scenes, we will nevertheless stick to the elementary language of class field theory in the main body of this~article. This~will turn out to be completely sufficient for our~purposes.
\end{ttt}

\section{Unramified primes}

\begin{propo}[{\rm Inert primes}{}]
\label{inert}
Let\/~$\frakl$
be a prime ideal
of\/~$K$
that is inert
in\/~$L/K$.
Write\/~$\ell := \#\calO_K/\frakl$,
denote by\/
$\frakL$
the unique prime
of\/~$L$
lying above\/
$\frakl$,
and assume~that

\begin{iii}
\item[ $\bullet$ ]
$a,b,c \in \calO_{L_\frakL}$,
$d \in \calO_{L_\frakL}^*$,
\item[ $\bullet$ ]
$(a/d \mod \frakL), (b/d \mod \frakL), (c/d \mod \frakL) \in \bbF_{\!\ell^3}$
are not contained
in\/~$\bbF_{\!\ell}$.
\end{iii}
Finally,~let\/
$S$
denote the surface~(\ref{zwei}).

\begin{abc}
\item[]
\begin{iii}
\item[{\rm a.i)} ]
If\/
$\ell > 3$
then\/
$S(K_\frakl) \neq \emptyset$.
\item[{\rm ii)} ]
If\/
$a \equiv b \pmod \frakl$
then\/
$S(K_\frakl) \neq \emptyset$.
\end{iii}
\item[{\rm b)} ]
For any\/
$(t_0\!:\!t_1\!:\!t_2\!:\!t_3) \in S(K_\frakl)$
such that\/
$t_0t_1 \neq 0$,
the quotient\/
$t_1/t_0 \in K_\frakl$
is in the image of the norm map
$N \colon L_\frakL \to K_\frakl$.
\end{abc}\medskip

\noindent
{\bf Proof.}
{\em
The~assumptions imply that
$a$,
$b$,
and~$c$
are
$\frakL$-adic
units,~too. Let~us write
$\overline{a} := (a \mod \frakL)$,
\ldots,
$\overline{d} := (d \mod \frakL) \in \bbF_{\!\ell^3}$.\smallskip

\noindent
a)
The reduction
$S_{\bbF_{\!\ell}}$
of
$S$
modulo~$\frakl$
is given by
$$T_0T_1T_2 = \N_{\bbF_{\!\ell^3}/\bbF_{\!\raisebox{-0.6mm}{\tiny $\ell$}}}(\overline{a} T_0 + \overline{b} T_1 + \overline{c} T_2 + \overline{d} T_3) \, .$$
It~is sufficient to show that
$S_{\bbF_{\!\ell}}$
admits a non-singular
$\bbF_{\!\ell}$-rational
point.\smallskip

\noindent
i)
We~claim that the singular locus
of~$S_{\overline\bbF_{\!\ell}}$
is zero-dimensional. To~show this, assume the~contrary. Then~there must exist a singular
point~$x$
on the
plane~$E := V(T_0)$.
Let~us write the equation
of~$S_{\overline\bbF_{\!\ell}}$
in the form
$T_0T_1T_2 = l_0l_1l_2$,
for
$l_0, l_1, l_2$
the three conjugates of the linear form
$\smash{\overline{a} T_0 + \overline{b} T_1 + \overline{c} T_2 + \overline{d} T_3}$.
Then,~by Lemma~\ref{geosing},
$x$
necessarily satisfies
$T_0 = T_1 = 0$
or
$T_0 = T_2 = 0$.
But~this would require
$\overline{c}/\overline{d} \in \bbF_{\!\ell}$
or~$\overline{b}/\overline{d} \in \bbF_{\!\ell}$
and therefore contradicts our~assumptions.

Hence,~the singular locus
of~$S_{\overline\bbF_{\!\ell}}$
is~finite. It~might happen that
$S_{\overline\bbF_{\!\ell}}$
is a cone over a smooth cubic curve, but then it certainly has a non-singular
\mbox{$\bbF_{\!\ell}$-rational}~point.
Otherwise,~it was shown by A.~Weil that
$\#S_{\overline\bbF_{\!\ell}}(\bbF_{\!\ell}) \geq \ell^2 - 2\ell + 1$
\cite[page~557]{Weil56}, cf.~\cite[Theorem~27.1 and Table~31.1]{Ma}. Further,~it is classically known that not more than four points may be singular~\cite[Corollary~9.2.3]{Dolgachev}.
As~$\ell^2 - 2\ell - 3 > 0$
for
$\ell > 3$,
this implies the~assertion.\smallskip

\noindent
ii)
We~have
$\overline{a} = \overline{b}$,
hence
$(1\!:\!(-1)\!:\!0\!:\!0) \in S_{\bbF_{\!\ell}}(\bbF_{\!\ell})$.
Lemma~\ref{geosing} shows that this is a non-singular~point.\smallskip

\noindent
b)
We~assume the coordinates of the point to be normalized such that
$t_0, \ldots, t_3 \in \calO_{K_\frakl}$
and at least one of them is a~unit. The~local extension
$L_\frakL/K_\frakl$
is unramified of degree~three. We~therefore have to show that
$\nu_\frakl(t_1/t_0)$
is divisible by~three.

Assume~the contrary.
If~$t_2 \neq 0$
then the equation of the surface ensures that
$3 | \nu_\frakl(t_0t_1t_2)$.
Thus,~the values
$\nu_\frakl(t_i)$,
for
$i = 0$,
$1$,
$2$,
must be mutually non-congruent
modulo~$3$.
Otherwise,~we know at least
$\nu_\frakl(t_0) \not\equiv \nu_\frakl(t_1) \pmod 3$
and~$t_2 = 0$.
There~are two~cases.\smallskip

\noindent
{\em First case.\/}
There is no unit among
$t_0, t_1, t_2$.

\noindent
Then~$t_3$
is a~unit. Since~we assume
$d$
to be a unit, too, we clearly have that
$at_0 + bt_1 + ct_2 + dt_3 \in \calO_{L_\frakL}^*$.
Hence,
$\N_{L_\frakL/K_\frakl}(at_0 + bt_1 + ct_2 + dt_3) \in \calO_{K_\frakl}^*$,
which, in view of
$t_0t_1t_2$
not being a unit, contradicts the equation of the surface.\smallskip

\noindent
{\em Second case.\/}
There is exactly one unit among
$t_0, t_1, t_2$.

\noindent
Without~restriction, assume that
$t_0$
is the~unit. Again,~we have that
$t_0t_1t_2$
is not a~unit. The~equation of the surface then requires that
$\N_{L_\frakL/K_\frakl}(at_0 + bt_1 + ct_2 + dt_3)$
must be a non-unit. To~ensure this, we need
$at_0 + bt_1 + ct_2 + dt_3 \not\in \calO_{L_\frakL}^*$,
which means nothing but
$$at_0 + dt_3 \equiv 0 \pmod \frakL \, .$$
But then
$a/d \equiv -t_3/t_0 \pmod \frakL$,
which is impossible since the right hand side
modulo~$\frakL$
is in
$\bbF_{\!\ell}$,
while the left hand side is~not.
}
\eop
\end{propo}

\begin{lem}
\label{geosing}
Let\/~$K$
be a field,
$l_0, l_1, l_2, l'_0, l'_1, l'_2 \in K[T_0,T_1,T_2,T_3]$
linear forms such that\/
$V(l_i) \neq V(l'_j)$
for all\/
$0 \leq i,j \leq 2$,
and\/
$S$
be the cubic surface, given
by\/~$l_0l_1l_2 = l'_0l'_1l'_2$.\smallskip

\noindent
Then,~every singular point
on\/~$S$
lying on the plane\/
$V(l_0)$
actually lies on at least two of the planes\/
$V(l_i)$
and two of the planes\/~$V(l'_j)$.\medskip

\noindent
{\bf Proof.}
{\em
Let~$x$
be a singular point
on~$S$
lying on the
plane~$V(l_0)$.
Then~$x \in V(l'_0l'_1l'_2)$.
Without~restriction, we may suppose that
$x \in V(l'_0)$.
Further,~the assumption ensures that, after a suitable change of coordinates, we may assume that
$l_0 = T_0$
and~$l'_0 = T_1$.
I.e.,~that
$S$
is given by the equation
$$F(T_0,T_1,T_2,T_3) := T_0l_1l_2 - T_1l'_1l'_2 = 0$$
and we consider a singular
point~$x = (0\!:\!0\!:\!t_2\!:\!t_3)$.
As~$\smash{\frac{\partial F}{\partial T_0}(x) = l_1l_2(x)}$
and
$\smash{\frac{\partial F}{\partial T_1}(x) = l'_1l'_2(x)}$,
the assertion~follows.
}
\eop
\end{lem}

\begin{lemo}[{\rm Split primes}{}]
\label{split}
Let\/~$\frakl$
be a prime
of\/~$K$
that is split
in\/~$L/K$,
$a,b,c,d \in L$
be arbitrary,
and\/~$S$
be the cubic surface, given by equation~(\ref{zwei}).

\begin{abc}
\item
Then\/~$S(K_\frakl) \neq \emptyset$.
\item
Every~non-zero element
of~$K_\frakl$
is a local norm for the extension\/
$L/K$
of global fields.
\end{abc}\medskip

\noindent
{\bf Proof.}
{\em
a)
The~scheme
$S_{K_\frakl}$
is defined by the equation
$$T_0T_1T_2 = \prod_{i=1}^3 [\sigma_i(a) T_0 + \sigma_i(b) T_1 + \sigma_i(c) T_2 + \sigma_i(d) T_3]$$
for
$\sigma_i\colon K \hookrightarrow L \!\otimes_K\! K_\frakl \to K_\frakl$
the three~homomorphisms. There~is a
\mbox{$K_\frakl$-rational}
line or plane, defined by
$T_0 = \sigma_1(b) T_1 + \sigma_1(c) T_2 + \sigma_1(d) T_3 = 0$.\smallskip

\noindent
b)
This~is a standard result from class field~theory.
}
\eop
\end{lemo}

\begin{remo}[{\rm The archimedean primes}{}]
\label{real}
\begin{iii}
\item
Let~$\sigma\colon K \hookrightarrow \bbR$
be a real~prime. Then,~for
$a$,
$b$,
$c$,
$d \in L$
arbitrary, we also have
$S_{\bbR,\sigma}(\bbR) \neq \emptyset$.
Further,~every non-zero real number is, with respect
to~$\sigma$,
a local norm for the extension\/
$L/K$
of global~fields.

Indeed,~as
$L/K$
is a cubic Galois extension, there are three real primes
$\smash{\sigma_i\colon L \to \bbR}$
extending~$\sigma$.
This~immediately implies the second~assertion. On~the other hand,
$S_{\bbR,\sigma}$
is given by the equation
$\smash{T_0T_1T_2 = \prod_{i=1}^3 [\sigma_i(a) T_0 + \sigma_i(b) T_1 + \sigma_i(c) T_2 + \sigma_i(d) T_3]}$.
Thus,~the same argument as above yields plenty of real~points.

\looseness-1
In~fact, it is known since the days of L.~Schl\"af\/li~\cite[pp.~114f.]{Schlafli} that a non-singular cubic surface
$\calS$
over~$\bbR$
always has real~points. A~nice geometric argument for this is as~follows. Start~with a
\mbox{$\bbC$-rational}
point
$x \in \calS(\bbC)$
not lying on any of the 27~lines.
Unless~$x$
is already the extension of a real point, there is a unique line connecting
$x$
with its complex conjugate
$\overline{x}$.
This~line meets
$\calS$
in a third point, which must be~real.
\item
For~$\sigma\colon K \hookrightarrow \bbC$
a complex prime and
$a$,
$b$,
$c$,
$d \in L$
arbitrary, we clearly have
$S_{\bbC,\sigma}(\bbC) \neq \emptyset$.
Further,~every non-zero complex number is a local norm with respect
to~$\sigma$.
\end{iii}
\end{remo}

\section{\mbox{Ramified \!primes--Reduction \!to \!a \!particular \!cone}}

\begin{ttt}
Local~class field theory shows that a local field with residue field
$\bbF_{\!\ell}$
for
$\ell \equiv 2 \pmod 3$
does not allow any ramified, cyclic cubic~extensions. Hence,~a cyclic cubic extension
$L/K$
may ramify only at
primes~$\frakl$
such that either
\mbox{$\#\calO_K/\frakl \equiv 1 \pmod 3$}
or
$\#\calO_K/\frakl$
is a power
of~$3$.
We~will consider the former case in~this~article.
\end{ttt}\pagebreak[3]

\begin{lem}
Let\/~$p \neq 3$
be a prime~number.

\begin{iii}
\item
Then~the equation
$$27T_0T_1T_2 = (T_0 + T_1 + T_2)^3$$
defines a cubic
curve\/~$C$
over\/~$\bbF_{\!p}$
with a node
at\/~$(1\!:\!1\!:\!1)$.
The~two tangent directions
at\/~$(1\!:\!1\!:\!1)$
are defined
over\/~$\bbF_{\!p}$
if and only
if\/~$p \equiv 1 \pmod 3$.
\item
Let\/~$e \geq 1$
be any~integer. Then,~for every\/
$\bbF_{\!p^e}$-rational
point\/
$(t_0\!:\!t_1\!:\!t_2)$
on\/~$C$,
at least one of the expressions\/
$t_1/t_0$,
$t_2/t_1$,
and\/
$t_0/t_2$
is properly defined and non-zero
in~$\bbF_{\!p^e}$.
Further,~these quotients evaluate solely to cubes
in~$\bbF_{\!p^e}^*$.
\end{iii}\medskip

\noindent
{\bf Proof.}
{\em
i)
It~is a standard calculation to show that
$C$
has a singular point
at~$(1\!:\!1\!:\!1)$
and no~others. The~tangent cone at this point is defined by a binary quadric of discriminant
$(-243) = -3 \!\cdot\! 9^2$.
Thus,~it splits if and only if
$(-3)$
is a quadratic residue
modulo~$p$.\smallskip

\noindent
ii)
The~first assertion simply says that
$(1\!:\!0\!:\!0), (0\!:\!1\!:\!0), (0\!:\!0\!:\!1) \not\in C$.
Further,~a calculation in any computer algebra system verifies that, in
$\bbZ[T_0, T_1, T_2]$,
the polynomial~expression
$$(T_0^2 + 2T_0T_1+T_1^2 + 5T_0T_2 - 4T_1T_2-5T_2^2)^3 + 729T_0(T_1-T_2)^3T_2^2$$
splits into two factors, one of which is
$27T_0T_1T_2 - (T_0 + T_1 + T_2)^3$.
Hence,~$T_0/T_2$
is the cube of a rational function
on~$C$.

Further,~for
$(t_0 \!:\! t_1 \!:\! t_2) \in C(\bbF_{\!p^e})$
with
$t_2 \neq 0$,
we see that
$t_0/t_2 \in \bbF_{\!p^e}$
is a cube, except possibly for the case
when~$t_1 = t_2$.
But~then the equation of the curve shows that
$t_0/t_2 = (\frac{t_0 + 2t_2}{3t_2})^3$.

Due~to symmetry, the same is true for
$t_1/t_0$
and~$t_2/t_1$.
This~is the~assertion.
}
\eop
\end{lem}

\begin{ex}
Consider the nodal cubic curve, defined by
$$T_0T_1T_2 + (T_0 + T_1 + T_2)^3 = 0$$
over~$\bbF_{\!7}$.
Besides~the node
at~$(1\!:\!1\!:\!1)$,
there are the six
\mbox{$\bbF_{\!7}$-rational}
points
$(1\!:\!(-1)\!:\!0), (1\!:\!0\!:\!(-1)), (0\!:\!1\!:\!(-1)), (1\!:\!1\!:\!(-1)), (1\!:\!(-1)\!:\!1)$,
and~$((-1)\!:\!1\!:\!1)$.
We~see explicitly that, for every
\mbox{$\bbF_{\!7}$-rational}
point, at least one of the expressions
$T_1/T_0$,
$T_2/T_1$,
$T_0/T_2$
is properly defined and non-zero
in~$\bbF_{\!7}$
and that all these quotients are cubes
in~$\bbF_{\!7}^*$.
\end{ex}

\begin{coro}
\label{nodal}
Let\/~$\ell$
be a prime power, but not a power of\/
$3$,
and\/
$\alpha \in \bbF_{\!\ell}^*$
arbitrary.

\begin{iii}
\item
Then~the equation\/
$T_0T_1T_2 - \frac1{27} (\alpha T_0 + T_1 + \frac1\alpha T_2)^3 = 0$
defines a nodal cubic
curve\/~$C'$
over\/~$\bbF_{\!\ell}$.
\item
For~every\/
\mbox{$\bbF_{\!\ell}$-rational}
point\/
$(t_0\!:\!t_1\!:\!t_2)$
on\/~$C'$,
at least one of the expressions\/
$t_1/t_0$,
$t_2/t_1$,
$t_0/t_2$
is properly defined and non-zero
in\/~$\bbF_{\!\ell}$.
Further,~these quotients evaluate only to elements in the coset
of\/~$\alpha$
modulo the cubic~residues.
\end{iii}\medskip

\noindent
{\bf Proof.}
{\em
i)
The curve
$C'$
is isomorphic
to~$C$
as it is obtained
from~$C$
by substituting
$\alpha T_0$
for~$T_0$
and
$\frac1\alpha T_2$
for~$T_2$.\smallskip

\noindent
This~implies, as well, the first assertion of~ii).
Further,~$T_1 / \alpha T_0 = \frac1\alpha T_1 / T_0$,
$\frac1\alpha T_2 / T_1$,
and
$\alpha T_0 / \frac1\alpha T_2 = \alpha^2 T_0 / T_2$
are cubes as soon as they are properly defined
in~$\bbF_{\!\ell}^*$.
}
\eop
\end{coro}

\begin{propo}[{\rm Ramified primes}{}]
\label{verzw}
Let\/~$\frakl$
be a prime ideal
of~$K$
that is ramified
in\/~$L/K$.
Suppose~that\/
$\ell := \#\calO_K/\frakl$
is not a power
of\/~$3$.
Denote~by\/
$\frakL$
the unique prime
of\/~$L$
lying
above\/~$\frakl$
and assume~that

\begin{iii}
\item[ $\bullet$ ]
$a \in \calO_{K_\frakL}, (a \mod \frakL) = \frac\alpha3$,
\item[ $\bullet$ ]
$b \in \calO_{K_\frakL}, (b \mod \frakL) = \frac13$,
\item[ $\bullet$ ]
$c \in \calO_{K_\frakL}, (c \mod \frakL) = \frac1{3\alpha}$,
\item[ $\bullet$ ]
$d \in \frakL \!\setminus\! \frakL^3$
\end{iii}
for some\/
$\alpha \in \bbF_{\!\ell}^*$. 
Finally,~let\/
$S$
denote the surface~(\ref{zwei}).

\begin{abc}
\item
Then~one has\/
$S(K_\frakl) \neq \emptyset$.
\item
Let\/
$(t_0\!:\!t_1\!:\!t_2\!:\!t_3) \in S(K_\frakl)$
be any~point. Then~not more than one
of\/~$t_0$,
$t_1$,
$t_2$
may~vanish.
\item[{\rm c.i)} ]
Suppose~that\/
$\alpha \in \bbF_{\!\ell}^*$
is a non-cube. Then~the following is~true.

Let\/
$(t_0\!:\!t_1\!:\!t_2\!:\!t_3) \in S(K_\frakl)$
be any~point. If,~for\/
$0 \leq i < j \leq 2$,
one has\/
$t_it_j \neq 0$
then the quotient\/
$t_j/t_i \in K_\frakl$
is\/ {\em not} in the image of the norm map\/
$N \colon L_\frakL \to K_\frakl$.
\item[{\rm ii)} ]
If\/
$\alpha \in \bbF_{\!\ell}^*$
is a cube then, for any\/
$(t_0\!:\!t_1\!:\!t_2\!:\!t_3) \in S(K_\frakl)$,
the quotients\/
$t_j/t_i \in K_\frakl$,
$0 \leq i < j \leq 2$,
are local norms, as soon as they are properly~defined.
\end{abc}\medskip

\noindent
{\bf Proof.}
{\em
Recall~that we automatically have
$\ell \equiv 1 \pmod 3$.\smallskip

\noindent
a)
The reduction
of~$S$
modulo~$\frakl$
is given by
$T_0T_1T_2 - \frac1{27} (\alpha T_0 + T_1 + \frac1\alpha T_2)^3 = 0$.
I.e.,~it is the cone over the nodal cubic
curve~$C'$,
studied in Corollary~\ref{nodal}.
$S_{\bbF_{\!\ell}}$~has exactly
$(\ell-1)\ell$
non-singular
\mbox{$\bbF_{\!\ell}$-rational}
points.\smallskip

\noindent
b)
We~assume the contrary and consider the case that
$t_1 = t_2 = 0$,
the others being~analogous. The~equation of the surface then requires
$at_0 + dt_3 = 0$.
As our assumptions imply
$a \neq 0$
and~$d \neq 0$,
we certainly have
$t_0 \neq 0$
and~$t_3 \neq 0$
and may write
$d/a = -t_0/t_3$.
Here,~the right hand side is an element
of~$K_\frakl$.
Hence,
$3 | \nu_\frakL(-t_0/t_3)$.
On~the other hand, since
$a$
is a unit,
$d/a \in \frakL \!\setminus\! \frakL^3$.
Thus,
$\nu_\frakL(d/a) = 1$
or~$2$,
which is a~contradiction.\smallskip

\noindent
c)
No~\mbox{$K_\frakl$-rational}
point
on~$S$
may reduce to the cusp
$(0\!:\!0\!:\!0\!:\!1) \in S_{\bbF_{\!\ell}}(\bbF_{\!\ell})$.
Indeed,~such a point could written in normalized form such that
$\nu_\frakl(t_0), \nu_\frakl(t_1), \nu_\frakl(t_2) \geq 1$
and
$t_3$
is a~unit. But~then
$\nu_\frakl(t_0t_1t_2) \geq 3$,
while
$$\nu_\frakl(\N_{L_\frakL/K_\frakl}(at_0 + bt_1 + ct_2 + dt_3)) = \nu_\frakL(at_0 + bt_1 + ct_2 + dt_3) = \nu_\frakL(dt_3) = \nu_\frakL(d) = 1 {\rm ~or~} 2 \, .$$

Consequently,~Corollary \ref{nodal}.ii) shows that at least one of the quotients
$t_j / t_i$
is properly defined and a unit
in~$\calO_{K_\frakl}$.
Further,~its residue
modulo~$\frakl$
is a cube if and only if
$\alpha$~is.
Hence,~$t_j/t_i \in K_\frakl$
is a local norm when
$\alpha$
is a cube and not a local norm,~otherwise. Remark~\ref{Brauer}.i) implies that the same is true for each of the three quotients, as soon as it is properly defined
in~$K_\frakl^*$.
}
\eop
\end{propo}

\section{The main result}
\label{sec4}

\begin{theo}
\label{main}
Let\/~$L/K$
be a cyclic cubic extension that is ramified at at least one~prime. Denote~by\/
$\frakp_1, \ldots, \frakp_r \subset \calO_K$
the ramified primes and
write\/~$\frakP_i$
for the unique prime lying
above\/~$\frakp_i$.
Suppose
that\/~$q_i := \#\calO_K/\frakp_i$
is not a power
of\/~$3$,
for
any\/~$i$.\smallskip

\noindent
Choose~a non-cube\/
$\alpha_1 \in \bbF_{\!q_1}^*$,
cubes
$\alpha_2 \in \bbF_{\!q_2}^*, \ldots, \alpha_r \in \bbF_{\!q_r}^*$,
and assume that\/
$a$,
$b$,
$c$,
$d \in \calO_L$
satisfy the following~conditions.

\begin{iii}
\item
$d$
splits as\/
$(d) = \frakP_1 \!\cdot\ldots\cdot\! \frakP_r \frakL_1 \!\cdot\ldots\cdot\! \frakL_s$,
where\/
$\N(\frakL_i)$
are prime
ideals\/~$\neq\! \frakp_1, \ldots, \frakp_r$.
I.e.,~$(d)$
does not contain any inert prime and contains\/
$\frakP_1, \ldots, \frakP_r$
exactly~once.
\item
$(a/d \mod \frakl\calO_L),\, (b/d \mod \frakl\calO_L),\, (c/d \mod \frakl\calO_L) \in \calO_L/\frakl\calO_L \setminus \calO_K/\frakl$
for every inert prime ideal\/
$\frakl$
in\/~$K$.
\item
$a \equiv b \pmod {\frakl\calO_L}$
for every inert prime\/
$\frakl$
of\/~$K$
such that\/
$\#\calO_K/\frakl = 2$
or\/~$3$.
\item
$(a \mod \frakP_i) = \frac{\alpha_i}3$,
$(b \mod \frakP_i) = \frac13$,
and\/
$(c \mod \frakP_i) = \frac1{3\alpha_i}$,
for\/~$i = 1, \ldots, r$.
\end{iii}
Finally,~let\/
$S$
denote the surface~(\ref{zwei}).
Then\/~$S(\bbA_K) \neq \emptyset$
but\/
$S(K) = \emptyset$.\medskip

\noindent
{\bf Proof.}
{\em
{\em First step.}
Preparations.

\noindent
Let~$\frakl \subset \calO_K$
be any prime~ideal.\smallskip

\noindent
{\em Case~1.}
$\frakl = \frakp_1, \ldots, \frakp_r$.

\noindent
Assumption~i) implies that,
for\/~$i = 1, \ldots, r$,
one has
$d \in \frakP_i \!\setminus\! \frakP_i^2 \subset \frakP_i \!\setminus\! \frakP_i^3$.
Together~with assumption~iv), we see that Proposition~\ref{verzw}~applies.\smallskip

\noindent
{\em Case~2.}
$\frakl$~is
an inert~prime.

\noindent
Then~$\frakL := \frakl\calO_L$
is the unique prime lying
above~$\frakl$.
By~assumption~i), we know that
$d \in \calO_{K_\frakL}^*$.
Thus,~we are exactly in the situation of Proposition~\ref{inert}.
Note~that we have
$S(K_\frakl) \neq \emptyset$
for primes
$\frakl$
with
$\#\calO_K/\frakl = 2,3$,
too, as then
$a \equiv b \pmod {\frakL}$.\smallskip

\noindent
{\em Case~3.}
$\frakl$~is
a split~prime.

\noindent
Then~Lemma~\ref{split}~applies.\smallskip

\noindent
{\em Second step.}
The existence of an adelic point.

\noindent
This~is equivalent to the existence of a real point
on~$S_{\bbR,\sigma}$
for every real prime
$\sigma\colon K \hookrightarrow \bbR$,
a complex point
on~$S_{\bbC,\sigma}$
for every complex prime
$\sigma\colon K \hookrightarrow \bbC$,
and~a
\mbox{$K_\frakl$-rational}
point for every prime ideal
$\frakl \subset \calO_K$.
We~know from Remark~\ref{real} that
$S_{\bbR,\sigma}$
and~$S_{\bbC,\sigma}$
admit real and complex points,~respectively. The~presence of
\mbox{$K_\frakl$-ra}\-tio\-nal
points for
every~$\frakl$
is guaranteed by Proposition~\ref{verzw}.a), Lemma~\ref{split}.a), and Proposition~\ref{inert}.a).\smallskip\pagebreak[3]

\noindent
{\em Third step.}
The non-existence of a
\mbox{$K$-rational}
point.

\noindent
Suppose~there would be a
{$K$-rational}
point
$(t_0\!:\!t_1\!:\!t_2\!:\!t_3) \in S(K)$.
Then,~by Proposition~\ref{verzw}.b),
among~$t_0, t_1, t_2$,
not more than one may~vanish. Choose
$0 \leq i < j \leq 2$
such that
$t_it_j \neq 0$
and consider the quotient
$t_j/t_i \in K$.

This~is a local norm at every archimedean prime and at every split prime by Remark~\ref{real} and Lemma~\ref{split}.b). It~is a local norm at every inert prime, too, in view of Lemma~\ref{inert}.b). Further,~the quotient
$t_j/t_i$
is a local norm at the ramified primes
$\frakp_2, \ldots, \frakp_r$,
as was shown in Proposition~\ref{verzw}.c). It~is, however, not a local norm at the
prime~$\frakp_1$.

The~proof is complete, since such a behaviour is in contradiction with global class field theory,~\cite[Chapter~VI, Corollary~5.7]{Neukirch} or~\cite[Theorem~5.1 together with~6.3]{Ta}.
}
\eop
\end{theo}

\begin{rem}
At~a first glance, condition~ii) seems to be hard to check as infinitely many primes are involed. Our~strategy for its verification is the~following.

\begin{iii}
\item[ $\bullet$ ]
Choose~elements
$z_1, z_2 \in \calO_L$
such that
$1, z_1, z_2$
are
\mbox{$K$-linearly}
independent.
Then~$\calO_L/\langle1,z_1,z_2\rangle$
is~finite. I.e.,~for an inert prime
$\frakl$
of~$K$
and
$\frakL$
the prime lying
above~$\frakl$,
one has
$\calO_{L_{\frakL}} = \langle1,z_1,z_2\rangle_{\calO_{K_\frakl}}$,
with only finitely many exceptions
$\frakl_1, \ldots, \frakl_s$.
\item[ $\bullet$ ]
For~$\frakl_1, \ldots, \frakl_s$,
check condition~ii)~directly.
\item[ $\bullet$ ]
Further,~write
$$\textstyle a/d = \frac1N(a_0 + a_1z_1 + a_2z_2), \quad
b/d = \frac1N(b_0 + b_1z_1 + b_2z_2), \quad
c/d = \frac1N(c_0 + c_1z_1 + c_2z_2)$$
for
$a_i, b_i, c_i \in \calO_K$
and
$N \in \calO_K$
an element that splits only into prime ideals split or ramified 
in~$L/K$
and, possibly,
into~$\frakl_1, \ldots, \frakl_s$.

Then~check that the ideals
$\gcd(a_1, a_2)$,
$\gcd(b_1, b_2)$,
$\gcd(c_1, c_2) \subset \calO_K$
do not contain any inert primes in their factorizations, except possibly
$\frakl_1, \ldots, \frakl_s$.
\end{iii}
\end{rem}

\begin{ex}
Let\/~$L := \bbQ(\zeta_7 + \zeta_7^{-1})$,
$z := \zeta_7 + \zeta_7^{-1} - 2$,
and
$S$
be the cubic surface
over~$\bbQ$,
given by equation~(\ref{zwei}), for
$$a := -1, \quad b := 5+6z^2, \quad c := 3+z^2, \quad d := z \, .$$
Then~$S$
violates the Hasse~principle.\medskip

\noindent
{\bf Proof.}
$L$~is
the unique cubic subfield of the cyclotomic extension
$\bbQ(\zeta_7)/\bbQ$.
It~is ramified only
at~$7$.
The~primes
$(\pm1 \mod 7)$
are split, the others are~inert. The~algebraic integer
$z$~is
chosen such that, for the unique prime
of~$L$
lying
above~$7$,
we have
$\frakP = (z)$.

We~work with the
\mbox{$\bbQ$-linearly}
independent elements
$1, z, z^2$,
which form a
\mbox{$\bbZ$-basis}
for~$\calO_K$.
Conditions~i) and~iii) of Theorem~\ref{main} are obviously satisfied. For~iv), note that,~in
$\bbF_{\!7}$,
one has
$\frac13 = 5$.
Further,~$\alpha = -\frac15 = 4$
is a non-cube and
$\frac54 = 3$.
Finally,
$a/d = \frac17(14 + 7z + z^2)$,
$b/d = \frac17(-70 + 7z - 5z^2)$,
and
$c/d = \frac17(-42 - 14z - 3z^2)$.
\eop
\end{ex}

\begin{exo}[{\rm continued}{}]
The~equation of
$S$
is, in explicit~form,
\begin{align*}
 & T_0^3 - 141T_0^2T_1 - 30T_0^2T_2 + 7T_0^2T_3 + 4863T_0T_1^2 + 2233T_0T_1T_2 - 532T_0T_1T_3 \\
 & {} + 251T_0T_2^2 - 119T_0T_2T_3 + 14T_0T_3^2 - 31499T_1^3 - 26286T_1^2T_2 + 6013T_1^2T_3 \\
 & {} - 6799T_1T_2^2 + 3157T_1T_2T_3 - 364T_1T_3^2 - 559T_2^3 + 392T_2^2T_3 - 91T_2T_3^2 + 7T_3^3 = 0 \, .
\end{align*}
$S$~has
bad reduction at
$2$,
$3$,
$7$,
$3739$,
and~$7589$.
While,
at~$7$,
the reduction is the cone over the nodal cubic curve, described in Corollary~\ref{nodal}, the reductions at the other bad primes only have one singular point~each.
$S_{\bbF_{\!2}}$
has a~binode. The~other three have a conical~singularity.

A~minimization algorithm yields a reembedding
of~$S$
as the surface, given by the equation
\begin{align*}
 & {-T_0^3} + 2T_0^2T_1 - T_0^2T_2 - 5T_0^2T_3 + T_0T_1^2 - T_0T_1T_2 + 7T_0T_1T_3 + 2T_0T_2^2 - 15T_0T_2T_3 \\
 & {}- 11T_0T_3^2 - T_1^3 - 2T_1^2T_2 + 9T_1^2T_3 + T_1T_2^2 + T_1T_3^2 + T_2^3 + T_2^2T_3 + 8T_2T_3^2 - T_3^3 = 0 \, .
\end{align*}
\end{exo}

\begin{coro}
\label{congr}
Let\/~$L/K$
be a cyclic cubic extension that is ramified at at least one prime, but unramified at all primes
dividing\/~$3$.
Choose~a split prime\/
$\frakl$
of\/~$K$
and let\/
$\smash{\widetilde{a}, \widetilde{b}, \widetilde{c}, \widetilde{d}}$
be four residue classes in\/
$[\calO_L / \frakl\calO_L]^* \cong (\bbF_{\!\ell}^*)^3$.\smallskip

\noindent
Then~there exists a cubic surface\/
$S$
that is a counterexample to the Hasse principle, of the~form
$$T_0T_1T_2 = \N_{L/K}(aT_0 + bT_1 + cT_2 + dT_3) \, ,$$
for\/
$a,b,c,d \in \calO_L$
such~that\/
$(a \mod \frakl\calO_L) = \widetilde{a}$,
\ldots,
$(d \mod \frakl\calO_L) = \widetilde{d}$.\medskip

\noindent
{\bf Proof.}
{\em
Denote~by
$\frakp_1, \ldots, \frakp_r \subset \calO_K$
the primes ramified
in~$L/K$.\smallskip

\noindent
{\em First step.}
$d$.

\noindent
Write~$\frakP_i$
for the unique prime ideal lying
above~$\frakp_i$,
for~$i = 1, \ldots, r$.
Further,~let
$\frakf \subseteq \calO_K$
be the conductor of the
extension~$L/K$
of global~fields.

According~to the Chinese remainder theorem, we may choose an element
$d' \in \calO_L$
such that
$(d' \mod \frakl\calO_L) = \widetilde{d}$
and
$d'$~is
a uniformizer for
each~$\frakP_i$.
We~find a partial~factorization
$$(d') = \frakP_1 \!\cdot\ldots\cdot\! \frakP_r \frakL_1 \!\cdot\ldots\cdot\! \frakL_s \!\cdot\! (m') \, ,$$
where
$\frakL_1, \ldots, \frakL_s$
are prime factors
in~$L$
of split primes
of~$K$
and
$m' \in \calO_K$
is an element that splits into a product of inert~primes. In~particular,
$m'$
is relatively prime
to~$\frakl$.

By~Sublemma~\ref{Cebo}, we may choose some
$m \equiv m' \pmod \frakl$,
$m \equiv 1 \pmod \frakf$,
such that
$(m)$
is a prime~ideal.
Indeed,~$\frakl$
and~$\frakf$
are relatively prime~\cite[Chapter~VI, Corollary~6.6]{Neukirch}, such that this is in fact a congruence condition
modulo~$\frakl\frakf$.
According to the decomposition law~\cite[Chapter~VI, Theorem~7.3]{Neukirch}, the congruence
$m \equiv 1 \pmod \frakf$
is clearly enough to ensure that
$(m)$
splits
in~$L/K$.

Thus,~put
$\smash{d := m \!\cdot\! \frac{d'}{m'}}$.
Then~$d \in \calO_L$
fulfills the congruence condition 
modulo~$\frakl\calO_L$
and assumption~i) of Theorem~\ref{main}.\smallskip

\noindent
{\em Second step.}
$c$.

\noindent
By~Lemma~\ref{congsys}, there is a solution
$c' \in \calO_L$
of the congruence~system
\begin{eqnarray*}
(c' \mod \frakl\calO_L) & = & \textstyle (\frac{\N(d)}d \mod \frakl\calO_L) \!\cdot\! \widetilde{c} \, , \\
(c' \mod \frakp_i\calO_L) & = & \textstyle \frac1{3\alpha_i} \!\cdot\! (\frac{\N(d)}d \mod \frakp_i\calO_L) \, , \quad i = 1,\ldots,r \, ,
\end{eqnarray*}
such that
$(c' \mod \frakp\calO_L) \not\in \calO_L/\frakp\calO_L \setminus \calO_K/\frakp$
for any prime
$\frakp$
of~$K$,
different from the ramified
ideals~$\frakp_j$.

Observe~here that
$\smash{\frac{\N(d)}d}$
is invertible
modulo~$\frakl\calO_L$
in view of the first~step. Further,
$\smash{\nu_{\frakP_i}(\frac{\N(d)}d) = 2}$.
Thus,~actually,
$\smash{(\frac{\N(d)}d \mod \frakp_i\calO_L) \in \frakP_i^2/\frakP_i^3}$,
which is a module over the residue~field.

Finally,~put
$\smash{c := c' \!\cdot\! \frac{d}{\N(d)}}$.
Then~$c \in \calO_L$
fulfills the congruence condition 
modulo~$\frakl\calO_L$
and the assumptions
on~$c$
made in Theorem~\ref{main}.ii) and~iv).\smallskip

\noindent
{\em Third step.}
$b$.

\noindent
Take~a solution
$b' \in \calO_L$
of the congruence~system
\begin{eqnarray*}
(b' \mod \frakl\calO_L) & = & \textstyle (\frac{\N(d)}d \mod \frakl\calO_L) \!\cdot\! \widetilde{b} \, , \\
(b' \mod \frakp_i\calO_L) & = & \textstyle \frac13 \!\cdot\! (\frac{\N(d)}d \mod \frakp_i\calO_L) \, , \quad i = 1,\ldots,r \, ,
\end{eqnarray*}
such that
$(b' \mod \frakp\calO_L) \not\in \calO_L/\frakp\calO_L \!\setminus\! \calO_K/\frakp$
for any prime ideal
$\frakp$
of~$K$,
different from
the~$\frakp_j$,
and put
$\smash{b := b' \!\cdot\! \frac{d}{\N(d)}}$.
Then~$b \in \calO_L$
fulfills the congruence condition 
modulo~$\frakl\calO_L$
and the assumptions
on~$b$
made in Theorem~\ref{main}.ii) and~iv).\smallskip

\noindent
{\em Fourth step.}
$a$.

\noindent
For~this, take a solution
$a' \in \calO_L$
of the slightly larger congruence~system
\begin{eqnarray*}
(a' \mod \frakl\calO_L) & = & \textstyle (\frac{\N(d)}d \mod \frakl\calO_L) \!\cdot\! \widetilde{b} \, , \\
(a' \mod \frakp_i\calO_L) & = & \textstyle \frac{\alpha_i}3 \!\cdot\! (\frac{\N(d)}d \mod \frakp_i\calO_L) \, , \quad i = 1,\ldots,r \, , \\
a' & \equiv & b' \pmod {\frakl_j} \, ,
\end{eqnarray*}
such that
$(a' \mod \frakp\calO_L) \not\in \calO_L/\frakp\calO_L \!\setminus\! \calO_K/\frakp$
for any unramified prime ideal
$\frakp$
of~$K$.
Here,~the
$\frakl_j \subset \calO_K$
run over all inert primes
of~$K$
such that
$\#\calO_K/\frakl_j = 2$
or~$3$.

To~complete the construction, put
$\smash{a := a' \!\cdot\! \frac{d}{\N(d)}}$.
Then~$a \in \calO_L$
fulfills the congruence condition 
modulo~$\frakl\calO_L$
and the assumptions
on~$a$
made in Theorem~\ref{main}.ii), iii), and~iv).
}
\eop
\end{coro}

\begin{lem}
\label{congsys}
Let\/~$K$
be a number field,
$L/K$
a cyclic cubic extension,
$\frakl_1, \ldots, \frakl_r \subset \calO_K$
distinct prime ideals, and\/
$a_1, \ldots, a_r \in \calO_K$.
Then,~for the congruence~system
\begin{eqnarray*}
a & \equiv & a_1 \pmod {\frakl_1\calO_L} \, ,\\
  & \ldots & \\
a & \equiv & a_r \pmod {\frakl_r\calO_L} \, .
\end{eqnarray*}
there~is a solution\/
$a \in \calO_L$
such that\/
$(a \mod \frakl\calO_L) \not\in \calO_L/\frakl\calO_L \!\setminus\! \calO_K/\frakl$
for any prime ideal\/
$\frakl \neq \frakl_1,\ldots,\frakl_r$
of\/~$K$,
unramified
in\/~$L/K$.\medskip

\noindent
{\bf Proof.}
{\em
Choose~elements
$z_1, z_2 \in \calO_L$
such that
$1, z_1, z_2$
are
\mbox{$K$-linearly}
independent.
Then~$\calO_L / \langle 1, z_1, z_2 \rangle$
is finite, such that
$\calO_L \!\otimes_{\calO_K}\! \calO_{K_\frakl} = \langle 1, z_1, z_2 \rangle_{\calO_{K_\frakl}}$
for all primes
$\frakl$
of~$K$,
with finitely many exceptions
$\frakp_1, \ldots, \frakp_s$.
Unless~these are among the
$\frakl_i$,
add congruence~conditions
$$a \equiv b_j \pmod {\frakp_j\calO_L} \, ,$$
for~$j \!=\! 1,\ldots,s$.
Choose
$b_1,\ldots,b_s \in \calO_L$
such that
$(b_j \mod \frakp_j\calO_L) \not\in \calO_L/\frakp_j\calO_L \setminus \calO_K/\frakp_j$,
whenever
$\frakp_j$
is unramified
in~$L/K$.

Since~$\calO_L$
is a Dedekind domain, the Chinese remainder theorem applies and we actually have only one congruence condition
$a \equiv A \pmod I$.
Take,~at first, any solution
$a' \in \calO_L$
of it and write it in the form
$$\textstyle a' = \frac1N(a_0 + a_1 z_1 + a_2 z_2) \, ,$$
for
$a_0, a_1, a_2 \in \calO_K$
and
$N \in \calO_K$
a product of the exceptional primes
$\frakp_1, \ldots, \frakp_s$.
The~congruence condition will clearly not be violated as long as we vary
$a_1$
and~$a_2$
in their respective residue classes
modulo~$NI$.

By~Sublemma~\ref{Cebo}, we may choose a representative
$\underline{a}_1$
such that
$(\underline{a}_1)$
is a product of some prime divisors
of~$NI$
and one further
prime~$\frakp'$.
In~addition, we may choose
$\underline{a}_2$
such that
$(\underline{a}_2)$
is a product of some
some prime divisors
of~$NI$
and a
prime~$\frakp'' \neq \frakp'$.
Then~$a := \frac1N(a_0 + \underline{a}_1 z_1 + \underline{a}_2 z_2) \in \calO_L$
solves the congruence system modulo the ideals
$\frakl_i$
and~$\frakp_j$.

Now let
$\frakl \subset \calO_K$,
$\frakl \neq \frakl_1,\ldots,\frakl_r$,
be any prime ideal, unramified
in~$L/K$.
If~$\frakl = \frakp_j$,
for
some~$j$,
then
$(a \mod \frakl\calO_L) \not\in \calO_L/\frakl\calO_L \setminus \calO_K/\frakl$
is fulfilled by~construction. But,~otherwise,
$\frakl$
is not a divisor of
of~$NI$.
Then~$\frac1N \in \calO_{K_\frakl}$
and
$\underline{a}_1$
and~$\underline{a}_2$
cannot be both divisible
by~$\frakl$.
This~implies the~assertion.
}
\eop
\end{lem}

\begin{subl}
\label{Cebo}
Let\/~$K$
be a number field,
$I \subset \calO_K$
an ideal,
and\/~$x \in \calO_K$.
Put\/~$\frakt := \lcm\limits_{n\in\bbN}(\gcd(I^n, (x)))$.\smallskip

\noindent
Then~there exist infinitely many pairwise non-associated elements\/
$y_1, y_2, \ldots \in \calO_K$
such that\/
$y_i \equiv x \pmod I$
and that each\/
$(y_i)$
factors into\/
$\frakt$
and a prime~ideal.\smallskip

\noindent
{\bf Proof.}
{\em
Write~$(x) = \frakx\frakt$.
Then~$\frakx$
is an ideal, relatively prime
to\/~$I$.
It~is known that there exist infinitely many prime ideals\/
$\frakp_i \subset \calO_K$
with the property~below.

There~exist some\/
$u_i, v_i \in\calO_K$,
$u_i \equiv v_i \equiv 1 \pmod I$
such~that
$$\frakp_i \!\cdot\! (u_i) = \frakx \!\cdot\! (v_i) \, .$$
Indeed,~the invertible ideals
of~$K$
modulo the principal ideals generated by elements from the residue class
$(1 \mod I)$
form an abelian group that is canonically isomorphic to the ray class group
$Cl_K^I \cong C_K/C_K^I$
of~$K$
\cite[Chapter~VI, Proposition~1.9]{Neukirch}. Thus,~the claim follows from the Cebotarev density theorem applied to the ray class field
$K^I\!/\!K$,
which has the Galois group
$\Gal(K^I\!/\!K) \cong Cl_K^I$.

Take~one of these prime~ideals.
Then~$\frakp_i\frakt \!\cdot\! (u_i) = \frakx\frakt \!\cdot\! (v_i) = (xv_i)$.
As~$\frakp_i\frakt \subset \calO_K$,
this shows that
$xv_i$
is divisible
by~$u_i$.
Put~$y_i := xv_i/u_i$.
Then~$(y_i) = \frakp_i\frakt$.
Further,
\mbox{$y_i \equiv x \pmod I$}.
}
\eop
\end{subl}

\begin{theo}
Let\/~$K$
be any number field,
$U_\reg \subset \Pb_K^{19}$
the open subset parametrizing non-singular cubic surfaces, and\/
$\calHC_K \subset U_\reg(K)$
be the set of all cubic surfaces
over\/~$K$
that are counterexamples to the Hasse~principle.\smallskip

\noindent
Then~the image of\/
$\calHC_K$
under Clebsch's invariant map
$$\smash{\C\colon U_\reg \longrightarrow \Pb(1,2,3,4,5)_K}$$
is Zariski~dense.\medskip

\noindent
{\bf Proof.}
{\em
Consider the family
$p\colon \calS \to \bbA^{12}_\bbC$
of cubic surfaces, given by the equation
$T_0T_1T_2 = \prod_{i=0}^2 \sum_{j=0}^3 a_{ij}T_j$.
Clebsch's~fundamental invariants, when applied to this family, define a rational~map
$\Cl\colon \bbA^{12}_\bbC \,\ratarrow\, \Pb(1,2,3,4,5)_\bbC$.
We~know that
$\Cl$
is~dominant. Indeed,~up to isomorphy, every non-singular cubic surface appears as a fiber of the
family~$p$~\cite[Corollary~9.3.4]{Dolgachev}.

For~a cyclic cubic
extension~$L/K$,
there is the similar family
$\smash{p^{(L)}\colon \calS^{(L)} \to \bbA^{12}_K}$,
given by
\begin{equation}
\label{NormEq}
T_0T_1T_2 = \smash{\prod_{i=0}^2 \left[ (a_0 + \sigma_i(z_1) a_1 + \sigma_i(z_2) a_2) T_0 + \cdots + (d_0 + \sigma_i(z_1) d_1 + \sigma_i(z_2) d_2)T_3 \right] \, ,}
\end{equation}
where
$(1,z_1,z_2)$
is a
\mbox{$K$-basis}
of~$L$
and
$\sigma_0, \sigma_1, \sigma_2 \in \Gal(L/K)$
denote the three~elements.

After~base extension
to~$\bbC$,
the family
$p^{(L)}$
becomes isomorphic
to~$p$.
As~the property of being dominant may be tested after extension of the base field, we see that the Clebsch invariant map
$\Cl^{(L)}\colon \bbA^{12}_K \,\ratarrow\, \Pb(1,2,3,4,5)_K$,
associated to the
family~$p^{(L)}$,
is dominant,~too.

Now~assume that
$\C(\calHC_K) \subset \Pb(1,2,3,4,5)_K$
were not Zariski~dense. Then, even more, the image
under~$\Cl^{(L)}$
of the counterexamples to the Hasse principle, contained in the
family~$p^{(L)}$,
has to be contained in a proper closed subset
$V \subset \Pb(1,2,3,4,5)_K$.
Since~$\Cl^{(L)}$
is dominant, this implies that
$\smash{\overline{(\Cl^{(L)})^{\raisebox{0.8mm}{\scriptsize -1}}(V)}}$
is a proper closed subset
of~$\bbA^{12}_K$.

In~other words, there exists a non-zero polynomial
$f \in K[A_j, \ldots, D_j]_{j=0,1,2}$
of a certain
degree~$d$
such that, for every counterexample to the Hasse principle of the form~(\ref{zwei}) with
$a = a_0 + a_1z_1 + a_2z_2$,
\ldots,
$d = d_0 + d_1z_1 + d_2z_2 \in \calO_L$,
one~has
$$f(a_0, a_1, a_2, b_0, b_1, b_2, c_0, c_1, c_2, d_0, d_1, d_2) = 0 \, .$$
Without~restriction, assume that the coefficients of
$f$
are algebraic~integers.

We~will show that this is in contradiction with our results~above. For~this, let us choose a particular
field~$L$
that it is ramified at at least one prime, but unramified at all primes lying
above~$3$.
Such~a choice is possible due to Lemma~\ref{extex}.

Further,~take a prime ideal
$\frakl \subset \calO_K$
that does not divide all the coefficients
of~$f$,
guarantees
$\calO_L \!\otimes_{\calO_K}\! \calO_{K_\frakl} = \langle 1, z_1, z_2 \rangle_{\calO_{K_\frakl}}$,
splits
in~$L$,
and is large enough to ensure
$(\ell-1)^{12} > d\ell^{11}$
for~$\ell := \#\calO_K/\frakl$.
The~existence is of such a prime follows from the decomposition law together with the Cebotarev density~theorem.

By~Corollary~\ref{congr}, we know that there are counterexamples to the Hasse principle of the form~(\ref{zwei}), with
$a,b,c,d \in \calO_L$
and
$(a \mod \frakl\calO_L)$,
\ldots,
$(d \mod \frakl\calO_L) \in (\calO_L/\frakl\calO_L)^*$
arbitrary. Consequently,
$$f(a_0, a_1, a_2, b_0, b_1, b_2, c_0, c_1, c_2, d_0, d_1, d_2) \equiv 0 \pmod \frakl \, ,$$
whenever
$a_0 + a_1z + a_2z^2$,
\ldots,
$d_0 + d_1z + d_2z^2 \in \calO_L$
are invertible
modulo~$\frakl\calO_L$.
This~shows that
$(f \mod \frakl)$
vanishes on at least
$(\ell-1)^{12}$
vectors
in~$\bbF_{\!\ell}^{12}$,
a contradiction to the lemma~below.
}
\eop
\end{theo}

\begin{lem}
Let\/~$\ell$
be a prime power and\/
$f \in \bbF_{\!\ell}[X_1,\ldots,X_n]$
a non-zero polynomial of
degree\/~$d$.\smallskip

\noindent
Then~the number\/
$N(\ell)$
of solutions of\/
$f(x_1,\ldots,x_n) = 0$
in\/~$\bbF_{\!\ell}^n$
satisfies the inequality\/
$N(\ell) \leq d\ell^{n-1}$.\medskip

\noindent
{\bf Proof.}
{\em
For~$\bbF_{\!\ell}$
the prime field, this is the lemma in \cite[Chapter~1, Paragraph~5.2]{BS}. For~the general case, the argument given there works equally~well.
}
\eop
\end{lem}

\begin{lem}
\label{extex}
Let\/~$K$
be a number field. Then~there exists a cyclic cubic extension\/
$L/K$
that is ramified at at least one prime
of\/~$K$,
but unramified at all primes
above\/~$3$.\medskip

\noindent
{\bf Proof.}
{\em
Probably~the easiest way to see this is as~follows.
Let~$p \equiv 1 \pmod 3$
be a prime number such that
$K/\bbQ$
is unramified
at~$p$.
Choose~$F$
to be the unique cubic subfield
of~$\bbQ(\zeta_p)/\bbQ$.
Then~$KF/\bbQ$
is ramified
at~$p$,
which shows
that~$L := KF \neq K$.
This~immediately implies
$\Gal(L/K) \cong \Gal(F/\bbQ) \cong \bbZ/3\bbZ$.
By~construction.
$L/K$
is unramified at all primes
above~$3$.
Further,~it must be ramified at some of the primes
of~$K$
lying
above~$p$.
}
\eop
\end{lem}

\begin{remo}[{\rm Variants}{}]
The~family constructed in Theorem~\ref{main} turned out to be sufficient to prove the main result, but it is certainly not unique in this~respect. At~least the following modifications are~possible.

\begin{iii}
\item
One~may allow that
$(d)$
contains some of the ideals
$\frakP_1$,
\ldots,
$\frakP_r$
twice instead of~once.
\item
The~congruence conditions required modulo the primes
$\smash{\frakP_1, \ldots, \frakP_r}$
could be~weakened. In~fact, the choice
of~$\alpha_i$
in Theorem~\ref{main} actually determines whether the local invariant
$(T_j/T_i, L_{\frakP_i}/K_{\frakp_i})$
is
$\frac13$,
$\frac23$,
or~$0$.
One~has to combine the conditions at the
$r$~ramified
primes in such a way that the sum is non-zero
in~$\frac13\bbZ/\bbZ$.
\end{iii}
\end{remo}

\appendix
\section{Other reduction types}

\begin{ttt}
One might ask for counterexamples to the Hasse principle of the form~(\ref{zwei}), having other reduction types at the ramified~primes. We~will show that this is impossible, at least
for~$\frakp$
not
dividing~$3$
and as long as we insist
on~$a,b,c,d \in \calO_{L_\frakP}$.
\end{ttt}

\begin{lem}
Let\/~$L$
be a cyclic cubic extension of the number
field\/~$K$
and\/
$\frakp$
be a prime
of\/~$K$
that is ramified
in\/~$L/K$.
Suppose~that\/
$q := \#\calO_K/\frakp$
is not a power
of\/~$3$
and
write\/~$\frakP$
for the unique prime
of\/~$L$
lying
above\/~$\frakp$.\smallskip

\noindent
Further,~let\/
$a,b,c,d \in L \cap \calO_{L_\frakP}$
and\/
$S$
be the cubic surface, given by~(\ref{zwei}). Assume that, for every\/
$(t_0\!:\!t_1\!:\!t_2\!:\!t_3) \in S(K_\frakp)$,
the quotients\/
$t_1/t_0$,
$t_2/t_1$,
$t_0/t_2$,
as soon as they are properly defined
in\/~$K_\frakp^*$,
are local~norms.\smallskip

\noindent
Then\/~$d \equiv 0 \pmod \frakP$
and\/~$27abc \equiv 1 \pmod \frakP$.\medskip

\noindent
{\bf Proof.}
{\em
The~assumption implies
$q \equiv 1 \pmod 3$.
Further,~the
reduction\/~$S_{\bbF_{\!q}}$
of the cubic
surface~$S$
at a ramified prime is given by
$$T_0T_1T_2 = (\overline{a}T_0 + \overline{b}T_1 + \overline{c}T_2 + \overline{d}T_3)^3 \, .$$
It~suffices to show that there are non-singular
\mbox{$\bbF_{\!q}$-rational}
points
on~$S_{\bbF_{\!q}}$
such that the quotients are non-cubes, except in the asserted~situation. There are three~cases.\smallskip

\noindent
{\em First case.}
$\overline{d} \neq 0$.

\noindent
Then, after a change of coordinates that does not involve
$T_0$
and~$T_1$,
$S_{\bbF_{\!q}}$
is the cubic surface, given by
$T_0T_1T_2 = T_3^3$.
It~is obvious that there are non-singular points such that
$t_1/t_0$
is a non-cube.\smallskip

\noindent
{\em Second case.}
$\overline{d} = 0$,
$\overline{a}, \overline{b}, \overline{c} \neq 0$.

\noindent
After~changing coordinates,
$S_{\bbF_{\!q}}$
is
$T_0T_1T_2 = A(T_0 + T_1 + T_2)^3$
for
$A = \overline{abc}$.
Assume~$A \neq \frac1{27}$,
as this is the claimed~exception.

Then~$S_{\bbF_{\!q}}$
is the cone over a non-singular
curve~$C$
of genus~one. The~triple cover
$T^3 = T_1/T_0$
is unramified and defines another curve
$\smash{\widetilde{C}}$
of genus~one. The assumption that
$t_1/t_0$
is always a cube leads to
$\smash{\#\widetilde{C}(\bbF_{\!q}) = 3 \!\cdot\! \#C(\bbF_{\!q})}$,
which contradicts Hasse's bound
for~$q \geq 16$.

Finally,~a systematic test shows that, for
$q = 4, 7, 13$
and
$A \neq 0, \frac1{27} \in \bbF_{\!q}$,
it does never happen that all the quotients are~cubes.\smallskip

\noindent
{\em Third case.}
$\overline{d} = 0$
and at least one of
$\overline{a}, \overline{b}, \overline{c}$
is~zero.

\noindent
Without~restriction, assume
$\overline{c} = 0$.
Then~the equation
of~$\smash{S_{\bbF_{\!q}}}$
simplifies to
$T_0T_1T_2 = (\overline{a}T_0 + \overline{b}T_1)^3$.
Considering~the partial derivative
by~$T_2$,
we see that, independently of what
$\overline{a}$
and~$\overline{b}$
are, every point such that
$t_0t_1 \neq 0$
is non-singular. We~may choose
$t_0, t_1 \in \bbF_{\!q}^*$
arbitrarily and find a point just by
calculating~$t_2$.
}
\eop
\end{lem}

\section{A geometric lemma}

\begin{lem}
\label{dim2}
The cubic surfaces having at least three Eckardt points are contained in a two-dimensional closed subset of the moduli scheme of non-singular cubic~surfaces.\medskip

\noindent
{\bf Proof.}
{\em
This~follows rather directly from the investigations undertaken by E.~Dardanelli and B.~van Geemen in~\cite{DG}. A~cubic surface over an algebraically closed base field may either allow a pentahedral form
$$a_0T_0^3 + \ldots + a_4T_4^3 = 0, \qquad T_0 + \ldots + T_4 = 0$$
or~not.

In~the first case, in order to have three Eckardt points, three of the five coefficients have to be equal to each other~\cite[2.2]{DG}. In~particular, the corresponding surfaces are contained in a two-dimensional subset of the moduli~scheme.

Otherwise,~the surface might be cyclic, ns1 or ns2~\cite[5.1--5.3]{DG}. Cyclic~surfaces form a one-dimensional subset, while ns2 surfaces form a two-dimensional subset in the moduli scheme~\cite[Theorem~6.6]{DG}. Finally,~the ns1 surfaces with at least three Eckardt points may be parametrized by a two-dimensional family~\cite[Proposition 5.7.(4) and~(5)]{DG}.
}
\eop
\end{lem}

\end{document}